\documentclass[11pt]{article}

	\usepackage{amsmath, amssymb, amscd, latexsym, epsfig}
\newcommand{\arrow}{\rightarrow}

\newcommand{\bb}{\mathbb}
\newcommand{\cx}{{\bb C}}

\newcommand{\half}{{\bb H}}

\newcommand{\reals}{{\bb R}}

\newcommand{\makefig}[3]{
	\begin{figure}[htbp]
        \refstepcounter{figure}
	\label{#2}
        \begin{center}
		~#3~\\
		\medskip
                {\sf Figure \thefigure.  #1}
        \end{center}
	\end{figure}
}

\thicklines
\newcommand{\maketab}[3]{
	\begin{figure}[htbp]
        \refstepcounter{figure}
	\label{#2}
        \begin{center}
		#3~\\
		\bigskip
                {\sf Table \thefigure.  #1}
        \end{center}
	\end{figure}
}

\renewcommand{\bold}[1]{\smallskip \noindent {\bf \boldmath #1 }\nopagebreak[4]}

\newcommand{\qed}{\nopagebreak[4]\hfill
\rule{2mm}{2.5mm} \bigskip \pagebreak[2]}

\newcommand{\mat}[1]{\left( \begin{smallmatrix} #1 \end{smallmatrix}\right) }
\newcommand{\Mat}[1]{\left(\begin{matrix} #1 \end{matrix}\right)}

\renewcommand{\tilde}{\widetilde}

\newcommand{\bs}{\backslash}
\newcommand{\bdry}{\partial}

\newcommand{\brackets}[1]{\langle #1 \rangle}

\newcommand{\closure}{\overline}

\newcommand{\isom}{\cong}

\newcommand{\mem}{\in}

\newcommand{\notmem}{\not\in}

\newcommand{\plusorminus}{\pm}

\newcommand{\superset}{\supset}

\newcommand{\AND}{\;\;\;\text{and}\;\;\;}

\newcommand{\GE}{\;\;\ge\;\;}

\newcommand{\EQ}{\;\;=\;\;}

\newcommand{\st}{\: : \:}         

\newcommand{\Dbar}{{\overline{D}}}
\newcommand{\Ebar}{{\overline{E}}}

\newcommand{\Ibar}{{\overline{I}}}

\newcommand{\Mbar}{{\overline{M}}}

\newcommand{\Ubar}{{\overline{U}}}

\newcommand{\Rhat}{{\widehat{\reals}}}

\newcommand{\chat}{{\widehat{\cx}}}

\newcommand{\rhat}{{\widehat{\reals}}}

\newcommand{\core}{\operatorname{core}}

\newcommand{\diam}{\operatorname{diam}}

\newcommand{\Fix}{\operatorname{Fix}}

\newcommand{\hull}{\operatorname{hull}}

\newcommand{\id}{\operatorname{id}}

\newcommand{\interior}{\operatorname{int}}

\newcommand{\Isom}{\operatorname{Isom}}

\newcommand{\Lie}{\operatorname{Lie}}

\newcommand{\M}{\operatorname{M}\!}

\renewcommand{\mod}{\operatorname{mod}}

\newcommand{\PSL}{\operatorname{PSL}}

\renewcommand{\Re}{\operatorname{Re}}

\newcommand{\SU}{\operatorname{SU}}

\newcommand{\T}{\operatorname{T}\!}

\newcommand{\cC}{{\cal C}}

\newcommand{\cP}{{\cal P}}

\newtheorem{theorem}{Theorem}[section]
\newtheorem{prop}[theorem]{Proposition}
\newtheorem{lemma}[theorem]{Lemma}
\newtheorem{cor}[theorem]{Corollary}

\makeatletter
\def\cleardoublepage{\clearpage\if@twoside \ifodd\c@page\else
    \thispagestyle{plain}\hbox{}\newpage\if@twocolumn\hbox{}\newpage\fi\fi\fi}

\def\ps@headings{\let\@mkboth\markboth
  \def\@oddfoot{}%
  \def\@evenfoot{}%
  \def\@evenhead{\small \sc\thepage\hfil\leftmark}
  \def\@oddhead{\small \sc \rightmark\hfil\thepage}
  \def\chaptermark##1{{
    \edef\@tempa{\ifnum \c@secnumdepth >\m@ne \@chapapp\ \thechapter. \fi}%
    \expandafter \markboth \expandafter{\@tempa ##1}{}}}%
  \def\schaptermark##1{\markboth {##1}{##1}}%
  \def\sectionmark##1{{
    \edef\@tempa{\ifnum \c@secnumdepth >\z@ \thesection. \fi}%
    \expandafter \markright \expandafter{\@tempa ##1}}}}

\def\thebibliography#1{\section*{References\@mkboth
 {References}{References}}\list
 {[\arabic{enumi}]}{\settowidth\labelwidth{[#1]}\leftmargin\labelwidth
 \advance\leftmargin\labelsep
 \usecounter{enumi}}
 \def\newblock{\hskip .11em plus .33em minus .07em}
 \sloppy\clubpenalty4000\widowpenalty4000
 \sfcode`\.=1000\relax}

\newif\if@restonecol
\def\theindex{\@restonecoltrue\if@twocolumn\@restonecolfalse\fi
\columnseprule \z@
\columnsep 35pt\twocolumn[\@makeschapterhead{Index}]
 \@mkboth{Index}{Index}\thispagestyle{plain}\parindent\z@
 \parskip\z@ plus .3pt\relax\let\item\@idxitem}
\def\@idxitem{\par\hangindent 40pt}

\def\endtheindex{\if@restonecol\onecolumn\else\clearpage\fi}

\def\footnoterule{\kern-3\p@ 
 \hrule width .4\columnwidth 
 \kern 2.6\p@} 
\@addtoreset{footnote}{chapter} 
\long\def\@makefntext#1{\parindent 1em\noindent 
 \hbox to 1.8em{\hss$^{\@thefnmark}$}#1}

\@addtoreset{equation}{section}

\renewcommand{\l@section}{\@dottedtocline{0}{1.5em}{2.3em}}
\renewcommand{\l@subsection}{\@dottedtocline{1}{3.8em}{3.2em}}
\renewcommand{\l@subsubsection}{\@dottedtocline{2}{7.0em}{4.1em}}

\makeatother

\newcommand{\Zk}{\RF_k M}

\newcommand{\F}{\operatorname{F}\!}
\newcommand{\RF}{\operatorname{RF}\!}

\newcommand{\gammat}{\tilde{\gamma}}
\newcommand{\ft}{\tilde{f}}

\newcommand{\zt}{\tilde{z}}
\newcommand{\bH}{\mathbb{H}}

\newcommand{\FM}{\F M}
\newcommand{\RFM}{\RF M}

\newcommand{\ba}{\backslash}

\newcommand{\Sier}{Sierpi\'nski~}

\title{ \vspace{-1in}
	{\bf Geodesic planes in the convex core of an acylindrical 3-manifold}
\vspace{.2in}}

\author{Curtis T. McMullen, Amir Mohammadi and Hee Oh}

\begin{document}

\maketitle

\begin{abstract}
Let $M$ be a convex cocompact, acylindrical hyperbolic 3-manifold of infinite volume,
and let $M^*$ denote the interior of the convex core of $M$.
In this paper we show that any geodesic plane in $M^*$ is either closed or dense.
We also show that only countably many planes are closed.  
These are the first rigidity theorems for planes in convex cocompact 3-manifolds
of infinite volume that depend only on the topology of $M$.
\end{abstract}

\tableofcontents

\vfill \footnoterule \smallskip
{\footnotesize \noindent
        Research supported in part by 
	the Alfred P. Sloan Foundation (A.M.) and the NSF.
	Typeset \today.
	}

\thispagestyle{empty}
\setcounter{page}{0}

\newpage

\section{Introduction}
\label{sec:intro}

In this paper we establish a new rigidity theorem for geodesic planes in 
acylindrical hyperbolic 3-manifolds.

\bold{Hyperbolic 3-manifolds.}
Let $M = \Gamma \bs \half^3$ be a complete, oriented hyperbolic $3$-manifold, presented as a quotient of hyperbolic space
by the action of a discrete group 
	$$\Gamma \subset  G =\Isom^+(\half^3).$$

Let $\Lambda \subset S^2 = \bdry \half^3$ denote the limit set of $\Gamma$,
and let $\Omega = S^2-\Lambda$ denote the domain of
discontinuity.  The {\em convex core} of $M$ is the smallest closed, convex subset of $M$ containing all closed geodesics;
equivalently,
\begin{displaymath}
	\core(M)=\Gamma\ba \hull(\Lambda)\subset M
\end{displaymath}
is the quotient of the convex hull of the limit set $\Lambda$ of $\Gamma$.
Let $M^*$ denote the interior of the convex core of $M$.

\bold{Geodesic planes in $M^*$.} 
Let 
\begin{displaymath}
        f : \half^2 \arrow M
\end{displaymath}
be a {\em geodesic plane}, i.e.\ a totally geodesic immersion of 
the hyperbolic plane into $M$.
We often identify a geodesic plane with its image, $P = f(\half^2)$.

By a geodesic plane $P^* \subset M^*$, we mean the nontrivial intersection
\begin{displaymath}
	P^* = P \cap M^* \neq \emptyset
\end{displaymath}
of a geodesic plane in $M$ with the interior of the convex core.
A plane $P^*$ in $M^*$ is always connected, and $P^*$ is closed in $M^*$ if and only if $P^*$ is
properly immersed in $M^*$ (\S\ref{sec:proper}).
 
\bold{Acylindrical manifolds and rigidity.}
In this work, we study geodesic planes in $M^*$ under the assumption that $M$ is
a convex cocompact, {\it acylindrical} hyperbolic $3$-manifold. 
The acylindrical condition is a topological one; it means that the compact Kleinian manifold
\begin{displaymath}
	\Mbar = \Gamma \bs (\half^3 \cup \Omega)
\end{displaymath}
has incompressible boundary, and every essential cylinder in $\Mbar$ is boundary parallel
(\S\ref{sec:proper}).
We will be primarily interested in the case where $M$ is a convex cocompact manifold of infinite volume.
Under this assumption, $M$ is acylindrical if and only if 
$\Lambda$ is a \Sier curve.\footnote{
A compact set $\Lambda \subset S^2$ is a {\em \Sier curve}
if $S^2-\Lambda = \bigcup D_i$
is a dense union of Jordan disks with disjoint closures,
and $\diam(D_i) \arrow 0$.
Any two \Sier curves are homeomorphic \cite{Whyburn:sier}.
}
	
Our main goal is to establish:

\begin{theorem} 
\label{thm:main}
Let $M$ be a convex cocompact, acylindrical, hyperbolic $3$-manifold.
Then any  geodesic plane $P^*$ in $M^*$ is either  closed or dense.
\end{theorem}

As a complement, we will show:

\begin{theorem}
\label{thm:countable}
There are only countably many closed geodesic planes $P^* \subset M^*$. 
\end{theorem}

We  also establish the following topological equidistribution result:
\begin{theorem}
\label{thm:dense}
If $P_i^*\subset M^*$ is an infinite sequence of distinct closed  geodesic planes, then
\begin{displaymath}
	\lim_{i\to \infty} P_i^*=M^*
\end{displaymath}
in the Hausdorff topology on closed subsets of $M^*$.
 \end{theorem}

\bold{Remarks.}
\begin{enumerate}
	\item
We do not know of any instance of Theorem \ref{thm:main}
where $P^*$ is closed in $M^*$ but $P$ is not closed in $M$.

{\em Added in proof.}   An example of such an {\em exotic plane}
in an acylindrical manifold has recently been constructed by Zhang.
In his example, the closure of $P$ is not even locally connected
near $\bdry M^*$
\cite{Zhang:exotic}.

Thus the rigidity of planes described in Theorem \ref{thm:main} does {\em not} 
extend beyond 
the convex core of $M$.  

	\item
In the special case where $M$ is compact (so $M=M^*$),
Theorem \ref{thm:main} is due independently to Shah and Ratner
(see \cite{Shah:immersions}, \cite{Ratner:topological}).

	\item
For a general convex cocompact manifold $M$, there can be uncountably
many distinct closed planes in $M^*$; see the end of \S\ref{sec:proper}.

	\item
Examples of acylindrical manifolds such that $M^*$
contains infinitely many closed geodesic planes 
are given in \cite[Cor.11.5]{McMullen:Mohammadi:Oh:sier} 

	\item 
The study of planes $P$ that do not meet $M^*$ can be reduced to the case where $M$ is
a quasifuchsian manifold.
This case can be analyzed via the bending lamination (cf. \S\ref{sec:near}).
\end{enumerate}

\bold{Comparison to the case of geodesic boundary.}
A convex cocompact hyperbolic $3$-manifold $M$ such that $\bdry \core(M)$ is totally geodesic is
automatically acylindrical.
For these {\em rigid} acylindrical manifolds, 
the results above were obtained in our previous work 
\cite{McMullen:Mohammadi:Oh:sier}.
While one would ultimately like to analyze planes in a large
class of geometrically finite groups, our previous results covered only
countably many examples (by Mostow rigidity).

The present paper makes a major step forward in this program,
by developing a new argument for unipotent recurrence
which works
{\em without} geodesic boundary,
which is robust enough to be invariant under quasi-isometry, and
which is powerful enough to apply to the class of all 
convex cocompact acylindrical manifolds.
The key insight is that one should work with a proper subset
of the renormalized frame bundle, defined in terms of thickness
of Cantor sets, where we show sufficient recurrence takes place
in the acylindrical case.

\makefig{Limit set of a cylindrical 3-manifold.}{fig:chars}{
\includegraphics[height=2.5in]{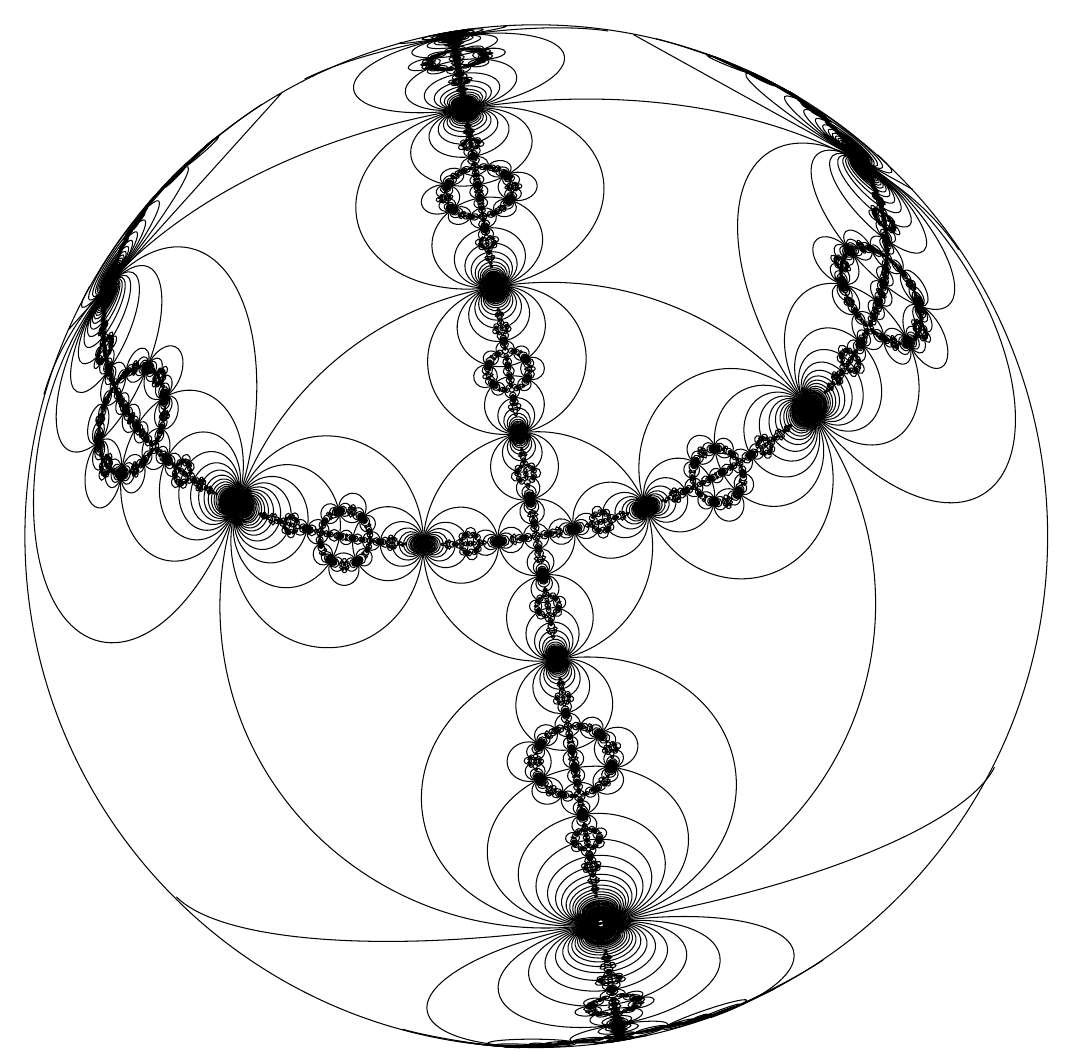}}

\bold{The cylindrical case.}
The acylindrical setting is also close to optimal, since
Theorem \ref{thm:main} is generally false
for cylindrical manifolds.

For example, consider a quasifuchsian group $\Gamma$ containing
a Fuchsian subgroup $\Gamma'$ of the second kind
with limit set $\Lambda' \subset S^1$.
Given $(a,b) \mem \Lambda' \times \Lambda'$, let $C_{ab}$ denote the unique circle orthogonal
to $S^1$ such that $C_{ab} \cap S^1 = \{a,b\}$.
It is possible to choose $\Gamma$ such that $C_{ab} \cap \Lambda = \{a,b\}$
for uncountably many $(a,b)$; and further, to arrange that the corresponding
hyperbolic planes $P \subset M$ and $P^* \subset M^*$ have wild closures,
violating Theorem \ref{thm:main}
(cf. \cite[App. A]{McMullen:Mohammadi:Oh:sier}).

The same type of example can be embedded in more complicated
3-manifolds with nontrivial characteristic submanifold; an example
is shown in Figure \ref{fig:chars}.  

\maketab{Notation for $G$ and some of its subgroups and 
homogeneous spaces.}{tab:Lie}{
\vspace{-.2in}
\begin{eqnarray*}
	G & = & \PSL_2(\cx) \isom \Isom^+(\half^3) \\
	H & = & \PSL_2(\reals) \isom \Isom^+(\half^2) \\
	K & = & \SU(2)/(\plusorminus I)\\
	A & = & \left\{ \mat{a&0 \\0 & a^{-1}} \st a > 0 \right\} \\
 	N & = & \left\{ n_s=\mat{1&s\\0&1} \st s \mem \cx \right\} \\
 	U & = & \{ n_s \st s \mem \reals \} \\
 	V & = & \{ n_s \st s \mem i\reals \}\\
	\F \half^3  & = & G = \{\text{the frame bundle of $\half^3$}\}\\
	\half^3 & = & G/K \\
	S^2 & = & G/AN = \bdry \half^3\\
	\cC & = & G/H = \{\text{the space of oriented circles $C \subset S^2$}\}\\
\end{eqnarray*}
\vspace{-.5in}
}

\bold{Homogeneous dynamics.}
Next we formulate a result in
the language of Lie groups and homogeneous spaces,
Theorem \ref{thm:homog},
that strengthens both Theorems \ref{thm:main} and \ref{thm:dense}. 

To set the stage, we have summarized our notation for
$G$ and its subgroups in Table \ref{tab:Lie}.
We have similarly summarized the spaces attached to an arbitrary
hyperbolic 3--manifold $M = \Gamma \bs \half^3$ in Table \ref{tab:M}.
(In the definition of $\cC^*$, a circle
$C \subset S^2$ {\em separates} $\Lambda$
if the limit set meets both components of $S^2-C$.)

\bold{Circles, frames and planes.}
Circles, frame and planes are closely related.  In fact, if $\cP$ denotes
the set of all (oriented) planes in $M$, then we have the natural identifications:
\begin{equation}
\label{eq:PCF}
	\cP = \Gamma \bs \cC = \FM/H .
\end{equation}
Indeed, all three spaces can be identified with $\Gamma \bs G / H$.
We will frequently use these identifications to go back and forth between
circles, frames and planes.

When $M^*$ is nonempty (equivalently, when $\Gamma$ is Zariski dense in $G$),
the spaces $\cC^*$ and $F^*$ correspond to the 
set of planes $\cP^*$ that meet $M^*$.  In other words, we have
\begin{equation}
\label{eq:star}
	\cP^*  = \Gamma \bs \cC^* = F^*/ H .
\end{equation}
To go from a circle to a plane, let $P$ be the image of $\hull(C) \subset \half^3$
under the covering map from $\half^3$ to $M$.  To go from
a frame $x \mem \FM$ to a plane, take the image of $xH$
under the natural projection $\FM \arrow M$.

When $\Lambda$ is connected and consists of more than one point
(e.g. when $M$ is acylindrical),
it is easy to see that:
\begin{displaymath}
	\closure{\cC^*} = \{C \mem \cC \st \text{$C$ meets $\Lambda$} \} .
\end{displaymath}
Thus the closures of the dense sets arising in Theorem 
\ref{thm:homog} below are quite explicit.

\maketab{Spaces associated to $M = \Gamma \bs \half^3$.}{tab:M}{
\vspace{-.2in}
\begin{eqnarray*}
	M & = & \Gamma \bs \half^3  =
		(\text{the quotient hyperbolic 3-manifold} )\\
	\Mbar & = & \Gamma \bs (\half^3 \cup \Omega) \\
	\core(M) & = & \Gamma \bs \hull(\Lambda)\\
	M^* & = & \interior(\core(M))\\
	\FM & = & \Gamma \bs G 
		= (\text{the frame bundle of $M$})\\
	F^* & = & \{x \mem \FM \st 
		\text{$x$ is tangent to a plane $P$ that meets $M^*$}\} \\
	\cC^* & = & \{C \mem \cC \st \text{$C$ separates $\Lambda$}\}\\
\end{eqnarray*}
\vspace{-.5in}
}

\bold{The closed or dense dichotomy.}
We can now state our main result from the perspective of 
homogeneous dynamics.

\begin{theorem}
\label{thm:homog} 
Let $M=\Gamma\ba \bH^3$ be a convex cocompact, acylindrical 3-manifold.
Then any $\Gamma$-invariant  subset of $\cC^*$ is either closed or  dense in $\cC^*$.
Equivalently, any $H$-invariant subset of $F^*$ is either closed or dense in $F^*$. 
\end{theorem}
(The equivalence is immediate from equation (\ref{eq:star}).)

This result sharpens Theorem \ref{thm:main} to give the 
following dichotomy on the level of the tangent bundles:
\begin{cor} 
\label{tangent}
The normal bundle to a geodesic plane $P^* \subset M^*$
is either closed or dense in the tangent bundle $\T M^*$.
\end{cor}

\bold{Beyond the acylindrical case.}
This paper also establishes several results that apply
outside the acylindrical setting.  For example,
Theorems \ref{thm:inc}, \ref{thm:zk}, \ref{thm:hull} and \ref{thm:case1} 
only require the assumption that $M$ has incompressible boundary.
In fact, the main argument pivots on a result 
relating Cantor sets and \Sier curves,
Theorem \ref{thm:inherit}, that involves no groups at all.

\bold{Discussion of the proofs.}
We conclude with a sketch of the proofs
of Theorems \ref{thm:main} through Theorem \ref{thm:homog}.

Let $M = \Gamma \bs \half^3$ be a convex cocompact acylindrical
3--manifold of infinite volume, with limit 
set $\Lambda$ and domain of discontinuity $\Omega$. 
The horocycle and geodesic flows on the frame bundle
$\FM = \Gamma \bs G$ are given by the right actions of $U$
and $A$ respectively.
The {\em renormalized frame bundle} of $M$ is the compact set defined by
\begin{equation}
\label{eq:RFM}
	 \RFM=\{x\in \FM : {xA}\text{ is bounded}\}.
\end{equation}

In \S\ref{sec:proper} we prove Theorem \ref{thm:countable}
by showing that the fundamental group of
any closed plane $P^* \subset M^*$ contains a free group
on two generators.
We also show that 
Theorems \ref{thm:main} and \ref{thm:dense} 
follow from Theorem \ref{thm:homog}.  The remaining
sections develop the proof of Theorem \ref{thm:homog}.

In \S\ref{sec:mod} we show that
$\Lambda$ is a \Sier curve of positive modulus.
This means there exists a $\delta > 0$ such that
the modulus of the annulus between any two 
components $D_1,D_2$ of $S^2-\Lambda$ satisfies
\begin{displaymath}
	\mod(S^2 - (\Dbar_1 \cup \Dbar_2)) \ge \delta > 0.
\end{displaymath}
We also show that if $\Lambda$ is a \Sier curve of
positive modulus, then there exists a $\delta > 0$
such that $C \cap \Lambda$ contains a Cantor set $K$
of modulus $\delta$, whenever $C$ separates $\Lambda$.
This means that for any disjoint components $I_1$ and $I_2$
of $C-K$, we have
\begin{displaymath}
	\mod(S^2 - (\Ibar_1 \cup \Ibar_2)) \ge \delta > 0.
\end{displaymath}
This result does not involve Kleinian groups and may be of interest
in its own right.

In \S\ref{sec:zk} we use this uniform bound on the modulus
of a Cantor set to construct a compact, $A$--invariant set
\begin{displaymath}
	\Zk \subset \RF M
\end{displaymath}
with good recurrence properties for the horocycle flow on $\FM$.
We also show that when $k$ is sufficiently large,
$\Zk$ meets every $H$--orbit in $F^*$.

The introduction of $\Zk$ is
one of the central innovations of this paper 
that allows us to handle acylindrical manifolds with quasifuchsian boundary.
When $M$ is a {\em rigid} acylindrical manifold,
$\Zk = \RF M$ for all $k$ sufficiently large, so in some
sense $\Zk$ is a substitute for the renormalized frame bundle.
For a more detailed discussion, see the end of \S\ref{sec:zk}.

In \S\ref{sec:hull} we shift our focus to the boundary of the convex core.
Using the theory of the bending lamination,
we give a precise description of $C \cap \Lambda$ in the case where $C$ comes
from a supporting hyperplane for the limit set.

In \S\S\ref{sec:near} and \ref{sec:far}, we formulate two
density theorems for hyperbolic 3-manifolds $M$ with incompressible boundary.
These results do not require that $M$ is acylindrical.
Each section gives a criterion for a sequence of circles $C_n \mem \cC^*$ to have
the property that $\bigcup \Gamma C_n$ is dense in $\cC^*$.

In \S\ref{sec:near} we show that density holds if $C_n \arrow C \notmem \cC^*$
and $\lim (C_n \cap \Lambda)$ is uncountable.  The proof relies on the
analysis of the convex hull given in \S\ref{sec:hull}.

In \S\ref{sec:far} we show that density holds if $C_n \arrow C \mem \cC^*$
and $C \notmem \bigcup \Gamma C_n$, provided $C \cap \Lambda$ contains a Cantor set of
positive modulus.  
The proof uses recurrence, minimal sets and homogeneous
dynamics on the frame bundle, and follows a similar argument
in \cite{McMullen:Mohammadi:Oh:sier}.  It also relies on the density result
of \S\ref{sec:near}.

When $M$ is acylindrical, the Cantor set condition is automatic by
\S \ref{sec:mod}.  Thus Theorem \ref{thm:homog} follows immediately from
the density theorem of \S\ref{sec:far}.

\bold{Question.}
We conclude by mentioning an open problem that goes
beyond the acylindrical case.
Let $P^* \subset M^*$ be a plane in a quasifuchsian manifold,
and suppose the corresponding circle satisfies $|C \cap \Lambda| > 2$.
Does it follow that $P^*$ is closed or dense in $M^*$?

\bold{Acknowledgements.}
We would like to thank Elon Lindenstrauss and Yair Minsky for useful discussions.

\section{Planes in acylindrical manifolds}
\label{sec:proper}

In this section we will prove Theorem \ref{thm:countable},
and show that our other main results, Theorems \ref{thm:main} and \ref{thm:dense},
follow from Theorem \ref{thm:homog} on the homogeneous dynamics of $H$ acting on $F^*$.

Let $M = \Gamma \bs \half^3$ be a convex cocompact hyperbolic 3-manifold.
We first describe how the topology of $\Mbar$ influences
the shape of planes in $M^*$.  Here are the two main results.

\begin{theorem}
\label{thm:inc}
If $\Mbar$ has incompressible boundary, then the fundamental group of any
closed plane $P^* \subset M^*$ is nontrivial.
\end{theorem}

\begin{theorem}
\label{thm:acy}
If $\Mbar$ is acylindrical, then the fundamental group of any
closed plane $P^* \subset M^*$ contains a free group on two generators.
\end{theorem}
The second result immediately implies Theorem \ref{thm:countable}, which we restate as follows:

\begin{cor}
\label{cor:countable}
If $\Mbar$ is acylindrical, then there are at most countably many
closed planes $P^* \subset M^*$.
\end{cor}

\bold{Proof.}
In this case $P^*$ corresponds to a circle $C$ whose stabilizer $\Gamma^C$
(as discussed below) is isomorphic to the fundamental group of $P^*$,
and contains a free group on two generators $\brackets{a,b}$.
Since $C$ is the unique circle containing the limit set of
$\brackets{a,b} \subset \Gamma$, and there are only countably many
possibilities for $(a,b)$, there are only countable possibilities for $P^*$.
\qed

In the remainder of this section, we first develop
general results about planes in 3--manifolds, and
prove Theorems \ref{thm:inc} and \ref{thm:acy}.
Then we derive Theorems \ref{thm:main} and \ref{thm:dense} from
Theorem \ref{thm:homog}.  Finally
we show by example that a cylindrical manifold
can have uncountably many closed planes $P^* \subset M^*$.

\bold{Topology of 3-manifolds.}
We begin with some topological definitions.

Let $D^2$ denote a closed 2-disk, and let $C^2 \isom S^1 \times [0,1]$ denote a closed cylinder.
Let $N$ be a compact 3-manifold with boundary.  
We say $N$ has {\em incompressible boundary} if every continuous map
\begin{displaymath}
	f : (D^2,\bdry D^2) \arrow (N,\bdry N)
\end{displaymath}
can be deformed, as a map of pairs, so its image lies in $\bdry N$.
(This property is automatic if $\bdry N = \emptyset$.)

Similarly, $N$ is {\em acylindrical} if it has incompressible boundary and every continuous map
\begin{displaymath}
	f : (C^2,\bdry C^2) \arrow (N,\bdry N),
\end{displaymath}
injective on $\pi_1$, can be deformed into $\bdry N$.
That is, every incompressible disk or cylinder in $N$ is boundary parallel.

When $N = \Mbar = \Gamma \bs (\half^3 \cup \Omega)$ is a compact Kleinian manifold,
these properties are visible on the sphere at infinity:  the limit set $\Lambda$ of $\Gamma$ is connected iff $\Mbar$ has incompressible
boundary, and $\Mbar$ is acylindrical iff $\Lambda$ is a \Sier curve or $\Lambda = S^2$.

For more on the topology of hyperbolic 3-manifolds, see e.g.\ 
\cite{Thurston:hype1}, \cite{Morgan:uniformization}, and \cite{Marden:book:HM}.

\bold{Topology of planes.}
Next we discuss the fundamental group of a plane $P \subset M$, and the corresponding plane $P^* \subset M^*$.
These definitions apply to an arbitrary hyperbolic 3-manifold.

For precision it is useful to think of a plane $P$ as being specified by an {\em oriented}
circle $C \subset S^2$, whose convex hull covers $P$.  More precisely, the plane attached to
$C$ is given by the map
\begin{displaymath}
	\ft : \hull(C) \isom \half^2 \subset \half^3 \arrow M = \Gamma \bs \half^3
\end{displaymath}
with image $\ft(\half^2) = P$.  The stabilizer of the circle $C$ in $G$ is a conjugate
$xHx^{-1}$ of $H = \PSL_2(\reals)$; hence its stabilizer in $\Gamma$ is given by
\begin{displaymath}
	\Gamma^C = \Gamma \cap x H x^{-1} .
\end{displaymath}

Let
\begin{displaymath}
	S = \Gamma^C \bs \hull(C) .
\end{displaymath}
Then the map $\ft$ descends to give an immersion 
\begin{displaymath}
	f : S \arrow M 
\end{displaymath}
with image $P$.
The immersion $f$ is generically injective if $P$ is orientable;
otherwise, it is generically two--to--one
(and there is an element in $\Gamma$ that reverses the orientation of $C$).

We refer to
\begin{displaymath}
	\pi_1(S) \isom \Gamma^C
\end{displaymath}
as the {\em fundamental group of $P$} (keeping in mind
caveats about orientability).

\bold{Planes in the convex core.}
Now suppose $P^* = P \cap M^*$ is nonempty.  In this case 
\begin{displaymath}
	S^* = f^{-1}(M^*)
\end{displaymath}
is a nonempty convex subsurface of $S$, with $\pi_1(S^*) = \pi_1(S)$.
The map
\begin{displaymath}
	f : S^* \arrow P^* \subset M^*
\end{displaymath}
presents $S^*$ as the (orientable) {\em normalization} of $P^*$,
i.e. as the smooth surface obtained by
resolving the self-intersections of $P^*$.  Similarly,
the frame bundle of $P$ with its branches separated is given by
\begin{displaymath}
	\F P = xH \subset \FM	
\end{displaymath}
for some $x \mem F^*$.  
(One should consistently orient $C$ and $P$ to define $\F P$.)

To elucidate the connections between these objects, we formulate:

\begin{prop}
\label{prop:discrete}
Let $M$ be an arbitrary hyperbolic 3--manifold.
Suppose $C \mem \cC^*$ and $x \mem F^*$ correspond to the
same plane $P^* \subset M^*$.
Then the following are equivalent:
\begin{enumerate}
	\item
$\Gamma C$ is closed in $\cC^*$.
	\item
The inclusion $\Gamma C \subset \cC^*$ is proper.
	\item
$xH$ is closed in $F^*$.
	\item
$P^*$ is closed in $M^*$.
	\item
The normalization map $f : S^* \arrow P^*$ is proper.
\end{enumerate}
\end{prop}
In (2) above, $\Gamma C$ is given the discrete topology.

\bold{Proof.}
If $\Gamma C$ is not discrete in $\cC^*$, then by homogeneity it is
perfect (it has no isolated points).  But a closed perfect set
is uncountable, so $\Gamma C $ is not closed.  Thus (1)
implies that $\Gamma C \subset \cC^*$ is closed and discrete, 
which implies (2); and clearly (2) implies (1).  
The remaining equivalences are similar, using
equation (\ref{eq:star}) to relate $\cP^*$, $\cC^*$ and $F^*$.
\qed

\bold{Compact deformations.}
In the context of proper mappings, the notion of a compact
deformation is also useful.

Let $f_0 : X \arrow Y$ be a continuous map.
We say $f_1 : X \arrow Y$ is a {\em compact
deformation} of $f_0$ if there is a continuous family
of maps $f_t : X \arrow Y$ interpolating between them,
defined for all $t\mem[0,1]$, and a compact set
$X_0 \subset X$ such that $f_t(x) = f_0(x)$ for all
$x \notmem X_0$.

Let $P^* \subset M^*$ be a hyperbolic plane with normalization
$f_0 : S^* \arrow M^*$.  We say $Q^* \subset M^*$ is a 
{\em compact deformation} of $P^*$ if it is the image of $S^*$
under a compact deformation $f_1$ of $f_0$.

\begin{theorem}
\label{thm:cc}
Let $M = \Gamma \bs \half^3$ be an arbitrary 3-manifold, and let
$K \subset M^*$ be a submanifold such that
the induced map 
\begin{displaymath}
	\pi_1(K) \arrow \pi_1(M)
\end{displaymath}
is surjective.
Then $K$ meets every geodesic plane $P^* \subset M^*$
and every compact deformation $Q^*$ of $P^*$.
\end{theorem}

\begin{cor}
\label{cor:K}
If $\pi_1(M)$ is finitely generated, then there is a compact
submanifold $K \subset M^*$ that meets every plane
$P^* \subset M^*$.
\end{cor}

\bold{Proof.}
Provided $M^*$ is nonempty, $\pi_1(M^*)$ is isomorphic to $\pi_1(M)$;
and since the latter group is finitely generated, 
there is a compact submanifold $K \subset M^*$ (say a neighborhood of
a bouquet of circles)
whose fundamental group surjects onto $\pi_1(M^*)$.
\qed

\bold{Proof of Theorem \ref{thm:cc}.}
We will use the fact that $S^0$ and $S^1$ can link in $S^2$.

Let $P^*$ be a plane in $M^*$, arising from a circle $C \subset S^2$
with an associated map $f : S \arrow P$ as above.
Since $P$ meets $M^*$, there are points in the limit set of $\Gamma$ on both sides of $C$.
Since the endpoints of closed geodesics are dense in $\Lambda \times \Lambda$ (cf. \cite{Eberlein:negative}),
we can find a hyperbolic element
$g \mem\Gamma$ such that its two fixed points 
\begin{displaymath}
	\Fix(g) = \{a_1,a_2\} \subset S^2
\end{displaymath}
are separated by $C$, and the convex hull of $\{a_1,a_2\}$ in $\half^3$ projects to a
closed geodesic $\delta \subset M$.  
Note that $\Fix(g) \isom S^0$ and $C\isom S^1$ are linked in $S^2$. 

Since $\pi_1(K)$ maps onto $\pi_1(M)$, the loop
$\delta$ is freely homotopic to
a loop $\gamma \subset K$.

Let $f_0 = f|S^*$.  Suppose $f_0 : S^* \arrow M^*$ has a compact deformation
$f_1$ with image $Q^*$ disjoint from $K$, and hence disjoint from
$\gamma$.  Extend this deformation trivially to the rest of $S$, to obtain
a compact deformation $f_1$ of the geodesic immersion $f : S \arrow P$.
Then $f_1(S)$ is disjoint from $\gamma$.  
Lifting $f_1$ to the universal cover of $S$,
we obtain a continuous map
\begin{displaymath}
	\ft_1 : \hull(C) \arrow \half^3
\end{displaymath}
that is a bounded distance from the identity map.
In particular, its image is a disk $D$ spanning $C$.

Similarly, a suitable lift of $\gamma$ gives a path
$\gammat \subset \half^3$, disjoint from $D$, that joins
$a_1$ to $a_2$.  This contradicts the fact that  
$C$ separates $a_1$ from $a_2$ in $S^2$.
\qed

We can now proceed to the:

\bold{Proof of Theorem \ref{thm:inc} (The incompressible case).}
For the beginning of the argument, we only use the fact that $\Mbar$ is compact
and $M^*$ is nonempty.  Using the nearest point projection, it is straightforward
to show that $\core(M)$ is homeomorphic to $\Mbar$.
Thus its interior $M^*$ deformation retracts onto a compact submanifold $K \subset M^*$,
homeomorphic to $\Mbar$, such that the inclusion is a homotopy equivalence;
in particular, $\pi_1(K) \isom \pi_1(\M^*)$.

Consider a closed plane
$P^* \subset M^*$, arising as the image of a proper map $f : S^* \arrow P^*$ as above.
We can also arrange that $K$ is transverse to $f$, so its preimage
\begin{displaymath}
	S_0 = f^{-1}(K) \subset S^*
\end{displaymath}
is a compact, smoothly bounded region in $S^*$.  (However $S_0$ need not be connected.)

We claim that, after changing $f$ by a compact deformation,
we can arrange that the inclusion of 
each component of $S_0$ into $S^*$ is injective
on $\pi_1$.  This is a standard argument in 3-dimensional topology.
If the inclusion is not injective on $\pi_1$, then there is a compact disk
$D \subset S^*$ with $D \cap S_0 = \bdry D$.
The map $f$ sends $(D,\bdry D)$ into $(M^*,K)$.
Since $K$ is a deformation retract of $M^*$,
$f|D$ can be deformed until it maps $D$ into $K$, while keeping $f|\bdry D$ fixed.
Then $D$ becomes part of $S_0$.  
This deformation is compact because $D$ is compact.
Since $\bdry S_0$ has only finitely many components,
only finitely many disks of this type arise, so after finitely many compact 
deformations of $f$, 
the inclusion $S_0 \subset S^*$ becomes injective on $\pi_1$.

Now we use the assumption that $K \isom \Mbar$ has incompressible boundary.
Suppose that $\pi_1(S^*)$ is trivial.
Then $\pi_1$ is trivial for each component
of $S_0$, and hence each component of $S_0$ is a disk.
By construction the deformed map $f$ restricts to give a map of pairs
\begin{displaymath}
	f : (S_0,\bdry S_0) \arrow (K,\bdry K).
\end{displaymath}
Since $K$ has incompressible boundary, we can further deform $f|S_0$ so it
sends the whole surface $S_0$ into $\bdry K$.  
Then the image $Q^*$ of $f$ gives a compact deformation of $P^*$
that is disjoint from $K^* = K-\bdry K$.
But $\pi_1(K^*)$ maps onto $\pi_1(M)$,
contradicting Theorem \ref{thm:cc}.
Thus $\pi_1(S^*)$ is nontrivial.
\qed

\bold{Proof of Theorem \ref{thm:acy} (The acylindrical case).}
The proof follows the same lines as the incompressible case.  If $\pi_1(S^*)$ does not
contain a free group on two generators, then $S^*$ is a disk or an annulus.
After a compact deformation, we can assume that the inclusion
$S_0 = f^{-1}(K) \subset S^*$ is injective on $\pi_1$.
Thus each component of $S_0$ is also a disk or an annulus.
Since $K$ is acylindrical, after a further compact deformation of $f$ 
we can arrange that $f(S_0) \subset \bdry K$, leading to a contradiction.
\qed

\bold{Rigidity of planes from homogeneous dynamics.}
Now suppose $M = \Gamma \bs \half^3$ is a convex cocompact, acylindrical 3-manifold.
Assume we know Theorem \ref{thm:homog}, which states that under this hypothesis: 
\begin{quote}
	{\em Any $\Gamma$--invariant set $E \subset \cC^*$ is closed or dense in $\cC^*$}.
\end{quote}
We can then prove the other two main results stated in the introduction.

\bold{Proof of Theorem \ref{thm:main}.}
Let $P^*$ be a geodesic plane in $M^*$, and let $E=\Gamma C$ be the corresponding set of circles.
Then by Theorem \ref{thm:homog}, $E$ is either closed or dense in $\cC^*$, and hence $P^*$ is either
closed or dense in $M^*$.
\qed

\bold{Proof of Theorem \ref{thm:dense}.}
Let $P_i^*$ be a sequence of distinct closed planes in $M^*$.
We wish to show that $\lim P_i^* = M^*$ in the Hausdorff
topology on closed subsets of $M^*$.
To see this, first pass to a subsequence so that $P_i^*$ converges to $Q^* \subset M^*$. 
It suffices to show that $Q^* = M^*$ for every such subsequence. 
Since each $P_i^*$ is nowhere dense, to show that $Q^* = M^*$ and complete the proof,
it suffices to show that $\bigcup P_i^*$ is dense in $M^*$.

Let $E_i \subset \cC^*$ be the $\Gamma$--orbit corresponding to
$P_i$, and let $E = \bigcup E_i$.
Since the planes $P_i$ are distinct, the sets $E_i$ are disjoint.
By Corollary \ref{cor:K}, there exists a compact set $K \subset M^*$ that meets every $P_i^*$, so there
exists a compact set $K' \subset \cC^*$ meeting every $E_i$.  
Thus we can choose $C_i \mem E_i \cap K'$ and 
pass to a subsequence such that 
\begin{displaymath}
	C_i \arrow C_\infty \mem K' \subset \cC^*
\end{displaymath}
and $C_\infty \notmem E$.  (If $C_\infty \mem E_i = \Gamma C_i$, just
drop that term from the sequence.)
Since $E$ is not closed in $\cC^*$, it is dense in $\cC^*$ by
Theorem \ref{thm:homog}.  Consequently $\bigcup P_i^*$ is dense in $M^*$,
as desired.
\qed

\bold{Example:  uncountably many geodesic cylinders.}
To conclude, we show that 
Theorem \ref{thm:acy} and Corollary \ref{cor:countable} do not hold for general convex cocompact manifolds with incompressible boundary.

In fact, in such a manifold one can have uncountably many distinct closed planes $P^* \subset M^*$, each with cyclic fundamental group.
For a concrete example of this phenomenon,
consider a closed geodesic $\gamma$ and the corresponding plane $P$ in the quasifuchsian manifold $M = M_\theta$
discussed in \cite[Cor. A.2]{McMullen:Mohammadi:Oh:sier}.  In this construction, $\gamma$ is a simple curve in the boundary of the convex core of $M$,
and $P \isom \gamma \times \reals$ is a hyperbolic cylinder properly embedded in $M$.
Consequently $P^* \subset M^*$ is a properly immersed cylinder in $M^*$.
By varying  the angle that $P$ meets the boundary of $\core(M_\theta)$ along $\gamma$, we obtain a continuous family
of properly immersed planes in $M^*$.

\section{Moduli of Cantor sets and \Sier curves}
\label{sec:mod}

The rest of the paper is devoted to the proof of Theorem \ref{thm:homog}.

In this section we define the modulus of 
a Cantor set $K \subset S^1$ (or in any circle $C \subset S^2$),
as well as the modulus of a \Sier curve $K \subset S^2$.
We then prove:

\begin{theorem}
\label{thm:robust}
Let $\Lambda$ be the limit set of $\Gamma$, where
$M=\Gamma \bs\half^3$ is a convex cocompact 
acylindrical 3--manifold of infinite volume.
Then there exists a $\delta>0$ such that:
\begin{enumerate}
	\item
$\Lambda$ is a \Sier curve of modulus $\delta$, and
	\item
$C \cap \Lambda$ contains a Cantor set of modulus $\delta$, whenever the circle $C\subset S^2$ separates $\Lambda$.
\end{enumerate}
\end{theorem}

\bold{The modulus of a \Sier curve.}
For background on conformal invariants and quasiconformal maps, see
\cite{Lehto:Virtanen:book}.

We begin with some definitions.
An {\em annulus} $A \subset S^2$ is an open region whose complement consists of
two components.  Provided neither component is a single point, $A$ is conformally
equivalent to a unique round annulus of the form
\begin{displaymath}
	A_R  = \{z \mem \cx \st 1 < |z| < R\},
\end{displaymath}
and its {\em modulus} is defined by
\begin{displaymath}
	\mod(A) = \frac{\log R}{2\pi} . 
\end{displaymath}
(More geometrically, $A$ is conformally equivalent to a Euclidean cylinder of radius $1$ and height $\mod(A)$.)
Since the modulus is a conformal invariant, we have
\begin{equation}
\label{eq:gmod}
	\mod(A) = \mod(g(A)) \;\;\forall g \mem G.
\end{equation}

Recall that a compact set $\Lambda \subset S^2$ is a {\em \Sier curve} if its complement
\begin{displaymath}
	S^2 - \Lambda = \bigcup D_i
\end{displaymath}
is a dense union of Jordan disks $D_i$ with disjoint closures, whose diameters
tend to zero.
We say $\Lambda$ has {\em modulus $\delta$} if 
\begin{displaymath}
	\inf_{i \neq j} \mod(S^2 - (\Dbar_i \cup \Dbar_j)) \ge \delta > 0.
\end{displaymath}

\bold{The modulus of an annulus $A \subset S^1$.}
Let $C \subset S^2$ be a circle and let $A \subset C$ be an `annulus on $C$',
meaning an open set such that $C-A = I_1 \cup I_2$ is the union of two disjoint intervals (circular arcs).
We extend the notion of modulus to this 1--dimensional situation by defining
\begin{displaymath}
	\mod(A,C) = \mod(S^2 - (I_1 \cup I_2)).
\end{displaymath}

Clearly $\mod(gA,gC) = \mod(A,C)$ for all $g \mem G$, and consequently
$\mod(A,C)$ depends only on the cross-ratio of the 4 endpoints of $A$.
The cross ratio is controlled by the lengths of the components
$A_1,A_2$ of $A$ and the components $I_1,I_2$ of $C-A$.
From this observation and monotonicity of the modulus \cite[I.6.6]{Lehto:Virtanen:book} it is easy to show:

\begin{prop}
\label{prop:xratio}
There are increasing continuous functions $\delta(t), \Delta(t) > 0$ such that 
\begin{displaymath}
	\delta(t) < \mod(A,C) < \Delta(t),
\end{displaymath}
where $t$ is the ratio of lengths
\begin{displaymath}
	t = \frac{\min(|A_1|,|A_2|)}{\min (|I_1|,|I_2|)}.
\end{displaymath}
\end{prop}
The same result holds with $t$ replaced by $d(\hull(I_1),\hull(I_2))$.

For later reference we recall the following result due to Teichm\"uller \cite[Ch II, Thm 1.1]{Lehto:Virtanen:book}:
\begin{prop}
\label{prop:Teich}
Let $I_1$ and $I_2$ be the two components of $C-A$.
Then
\begin{displaymath}
	\mod(B) \le \mod(A,C) 
\end{displaymath}
for any annulus $B \subset S^2$ separating 
the endpoints of $I_1$ from those of $I_2$.
\end{prop}

\bold{The modulus of a Cantor set.}
Let $K \subset C \subset S^2$ be a compact subset of a circle, such that its complement
\begin{displaymath}
	C - K = \bigcup I_i
\end{displaymath}
is a union of open intervals with disjoint closures.  
Note that $C$ is uniquely determined by $K$ (and we allow $K=C$).
We say $K$ has {\em modulus $\delta$} if we have
\begin{equation}
\label{eq:modK}
	\inf_{i \neq j} \mod(A_{ij},C) \ge \delta > 0,
\end{equation}
where $A_{ij} = C - \closure{I_i \cup I_j}$.
We will be primarily interested in the case where $K$ is a {\em Cantor set},
meaning $\bigcup I_i$ is dense in $C$.

\bold{Slices.}
Next we show that circular slices of a \Sier curve inherit positivity of the modulus.
This argument makes no reference to 3--manifolds.

\begin{theorem}
\label{thm:inherit}
Let $\Lambda \subset S^2$ be a \Sier curve of modulus $\delta > 0$.
Then there exists a $\delta'>0$ such that $C \cap \Lambda$ contains
a Cantor set $K$ of modulus $\delta'$ 
whenever $C$ is a circle separating $\Lambda$.
\end{theorem}

\makefig{A circle $C$ and some components $D_i$ of $S^2-\Lambda$.}{fig:blobs}{
\includegraphics[height=1.8in]{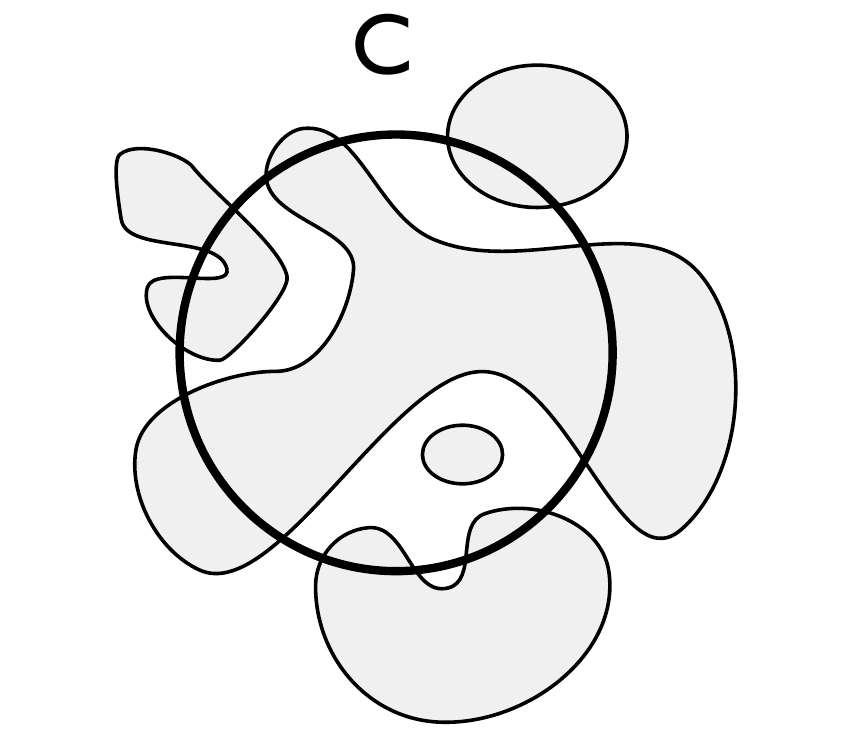}}

\bold{Proof.}  
Let $S^2 - \Lambda = \bigcup D_i$ express the complement of $\Lambda$ as a union of disjoint disks.
Each disk $D_i$ meets the circle $C$ in a collection of disjoint
open intervals (see Figure \ref{fig:blobs}).  The proof will be
based on a study of the interaction of intervals from
different components.

Let $U = C -\Lambda = \bigcup U_i$, where
\begin{displaymath}
	U_i = C \cap D_i .
\end{displaymath}
Note that distinct $U_i$ have disjoint closures,
and $\diam U_i \arrow 0$,
since these two properties hold for the disks $D_i$.  
The open set $U_i$ may be empty.

We may assume $U$ is dense in $C$, since otherwise we can just choose
a suitable Cantor set $K \subset C-\Ubar$.  On the other hand, no $U_i$ is dense in $C$;
if it were, we would have $C \subset \Dbar_i$, contrary to our assumption
that $C$ separates $\Lambda$.  It follows that $U_i$ is nonempty
for infinitely many values of $i$.

Let us say an open interval $I = (a,b) \subset C$, with distinct
endpoints, is a {\em bridge of type $i$} if $a,b \mem \bdry U_i$.
Note that an ascending union of bridges of type $i$ is 
again a bridge of type $i$, provided its endpoints are distinct.

Our goal is to construct a sequence of disjoint bridges
$I_1, I_2, I_3\ldots$ such that $|I_1| \ge |I_2| \ge \cdots$ and $K = C -\bigcup I_i$ 
is a Cantor set of modulus $\delta'$.

To start the construction, choose any bridge $I_1 \subset C$.
After changing coordinates by a M\"obius transformation $g\mem G^C$, we can assume
that $I_1$ fills at least half the circle; i.e.  $|I_1| > |C|/2$.
This will ensure that $|I_1| \ge |I_k|$ for all $k>1$.

Next, let $I_2$ be a bridge of maximal length among all those which are disjoint from $I_1$
and of a different type from $I_1$.
Such a bridge exists because $\diam(U_i) \arrow 0$, so only finitely many types of bridges are competing
to be $I_2$.  To complete the initial step, enlarge $I_1$ to a maximal interval of the same type,
disjoint from $I_2$.  

Proceeding inductively, let $I_{k+1}\subset C$ be a bridge
of maximum length among all bridges disjoint from $I_1,\ldots, I_k$.
Since $I_1$ is a maximal bridge of its type among those disjoint from $I_2$, and vice--versa,
the intervals $(I_1,I_2,I_k)$ are of 3 distinct types, for all $k \ge 3$.  Consequently
$|I_2| \ge |I_k|$ for all $k > 2$.

Note that the bridges so constructed have disjoint closures.
Indeed, if $I_i$ and $I_j$ were to have an endpoint $a$ in common,
with $i<j$, then $I_i \cup \{a\} \cup I_j$ would be a longer
interval of the same type as $I_i$, contradicting to stage $i$ of the construction.

Since $U$ is dense in $C$, it follows that at any finite stage there
is a bridge disjoint from all those chosen so far, and thus the
inductive construction continues indefinitely.  
By construction, we have
\begin{displaymath}
	|I_1| \ge |I_2| \ge |I_3| \cdots
\end{displaymath}
and by disjointness, $|I_k| \arrow 0$.  
Moreover, $\bigcup I_k$ is dense in $C$.
Otherwise, by density of $U$, we would be able to find a bridge $J$ disjoint from 
all $I_k$, and longer than $I_k$ for all $k$ sufficiently large,
contradicting the construction of $I_k$.

Let $K = C - \bigcup_1^\infty I_k$.
Since the intervals $I_k$ have disjoint closures,
and their union is dense in $C$,
$K$ is a Cantor set.  We have $K \subset \Lambda$ since
$\bdry I_k \subset \Lambda$ for all $k$.

Now consider any two indices $i < j$.
Let 
\begin{displaymath}
	A = C - (\Ibar_i \cup \Ibar_j) = A_1 \cup A_2,
\end{displaymath}
where the open intervals $A_1$ and $A_2$ are disjoint.
If the bridges $I_i$ and $I_j$ have types $s \neq t$ respectively,
then the annulus
\begin{displaymath}
	B = S^2 - (\Dbar_s \cup \Dbar_t)
\end{displaymath}
separates $\bdry I_i$ from $\bdry I_j$, and hence
\begin{displaymath}
	\mod(A,C) \ge \mod(B) \ge \delta > 0
\end{displaymath}
by Proposition \ref{prop:Teich}.

On the other hand, if $I_i$ and $I_j$ have the same type $s$,
then $i,j > 2$, and there must
be a bridge $I_k$, $k<i$, such that $I_1 \cup I_k$
separates $I_i$ from $I_j$.  Otherwise, we could
have combined $I_i$ and $I_j$ to obtain
a longer bridge at step $i$.

It follows that 
\begin{displaymath}
	t \EQ \frac{\min(|A_1|,|A_2|)}{\min (|I_i|,|I_j|)}
	\GE \frac{\min(|I_1|,|I_k|)}{\min (|I_i|,|I_j|)}
	\EQ \frac{|I_k|}{|I_j|}
	\GE 1,
\end{displaymath}
since $k < i < j$.
By Proposition \ref{prop:xratio}, this implies that 
\begin{displaymath}
	\mod(A,C) > \delta_0 > 0
\end{displaymath}
where $\delta_0$ is a universal constant.
Thus the Theorem holds with $\delta' = \min(\delta_0,\delta)$.
\qed

\bold{Limit sets.}
We can now complete the proof of Theorem \ref{thm:robust}.  

\begin{theorem}
\label{thm:Siermod}
Let $M = \Gamma \bs \half^3$ be a convex cocompact acylindrical 3--manifold of infinite volume.
Then its limit set $\Lambda$ is a \Sier curve of modulus $\delta$ for some $\delta>0$.
\end{theorem}

\bold{Proof.}
First suppose that every component of $\Omega = S^2-\Lambda = \bigcup D_i$ is a round disk, i.e.\
suppose that $M$ is a rigid acylindrical manifold.
By compactness, there exists an $L>0$ such that the hyperbolic length of
any geodesic arc $\gamma \subset \core(M)$ orthogonal to the boundary at its endpoints is greater
than $L$.  Consequently $d_{ij} = d(\hull(D_i),\hull(D_j)) \ge L$ for any $i \neq j$.
Since the modulus of $S^2-(\Dbar_i \cup \Dbar_j)$ is given by $d_{ij}/(2\pi)$,
$\Lambda$ is a \Sier curve of modulus $\delta = L/(2\pi) > 0$.

To treat the general case, recall that for any convex cocompact acylindrical manifold $M$,
there exists a rigid acylindrical manifold $M' = \Gamma' \bs \half^3$
such that $\Gamma'$ is $K$--quasiconformally conjugate to $\Gamma$.
Since a $K$--quasiconformal map distorts the modulus of an annulus by at most
a factor of $K$, and the limit set $\Lambda'$ of $\Gamma'$ is a \Sier curve with modulus $\delta' >0$,
$\Lambda$ itself is a \Sier curve of modulus $\delta = \delta'/K > 0$.
\qed

\bold{Proof of Theorem \ref{thm:robust}.}
Combine Theorems \ref{thm:inherit} and \ref{thm:Siermod}.
\qed

\section{Recurrence of horocycles}
\label{sec:zk}

Let $M = \Gamma \bs \half^3$ be an arbitrary 3--manifold.
In this section we will define, for each $k>1$, a closed, $A$--invariant set
\begin{displaymath}
	\Zk \subset \RF M
\end{displaymath}
consisting of points with good recurrence properties under the horocycle flow generated by $U$
(for terminology see Tables \ref{tab:Lie} and \ref{tab:M}).
We will then show:

\begin{theorem}
\label{thm:zk}
Let $M = \Gamma \bs \half^3$ be a convex cocompact acylindrical 3--manifold.
We then have
\begin{displaymath}
	F^* \subset (\Zk) H 
\end{displaymath}
for all $k$ sufficiently large.
More precisely, every plane $P^* \subset M^*$ is tangent to a frame in $\Zk$.
\end{theorem}
We conclude by comparing the general result above to results
that hold only when $\bdry M^*$ is totally geodesic.

We remark that $(\Zk)H$ is usually not closed, even when $M$ is acylindrical,
because there can be circles $C \mem \closure{\cC^*}$ such that
$|C \cap \Lambda| = 1$.

\bold{Thick sets.}
We begin by defining $\Zk$.  
Let us say a closed set $T \subset \reals$ is {\em $k$-thick} if
\begin{displaymath}
        [1,k] \cdot |T| = [0,\infty) .
\end{displaymath}
In other words, given $x \ge 0$ there exists a $t \mem T$ with
$|t| \mem [x,kx]$.
Note that if $T$ is $k$--thick, so is $\lambda T$ for all $\lambda \mem \reals^*$.

If the translate $T-x$ is $k$--thick for every $x \mem T$, we say $T$ is {\em globally} $k$--thick.
A set $K \subset U$ is (globally) $k$-thick if its image under an isomorphism $U \isom \reals$
is (globally) $k$--thick.

\bold{Unipotent recurrence.}
For $x \mem \RFM$, the unipotent orbit $xU$ almost never
remains in $\RFM$.   
Provided, however, there is a {\em thick set} $K \subset U$
such that $xK \subset \RFM$, we have sufficient recurrence to carry through
many arguments that would be automatic if $xU$ were bounded.
The key point is to combine thickness with the polynomial behavior of
unipotent flows.
This theme is developed in detail in \cite[\S 8]{McMullen:Mohammadi:Oh:sier},
and it motivates the definition of $\Zk$ below.

Let
\begin{equation}
\label{eq:Uz}
	U(z) = \{u \mem U \st zu \mem \RFM\}
\end{equation}
denote the return times of $z \mem \FM$ to the renormalized frame bundle
under the horocycle flow.
We define $\Zk$ for each $k>1$ by
\begin{displaymath}
	\Zk = \left\{z
		\mem \RFM \st
		\begin{array}{c}
		\text{there exists a globally $k$--thick}\\
		\text{set $K$ with $0 \mem K \subset U(z)$}
		\end{array}
		\right\}  .
\end{displaymath}
Let
\begin{displaymath}
	U(z,k) = \{u \mem U \st zu \mem \Zk\} .
\end{displaymath}

\begin{prop}
\label{prop:zkcompact}
Suppose the convex core of $M$ is compact.
Then for any $k>1$, the set $\Zk$ is a compact, $A$--invariant subset of $\RF M$.
Moreover, $U(z,k)$ is $k$--thick for each $z \mem \Zk$.
\end{prop}

\bold{Proof.}
Using compactness of $\RF M$, it is easily verified that if $z_n \arrow z$ in $\FM$ then
$\limsup U(z_n) \subset U(z)$.
One can also check that if $K_n \subset U$ is a sequence of globally
$k$--thick sets with $0 \mem K_n$, then $\limsup K_n$ is also globally $k$--thick.  
Consequently $\Zk \subset \RFM$ is closed, and hence compact.

Since $U(za)$ is a rescaling of $U(z)$ for any $a \mem A$, and the notion of thickness is
scale--invariant, $\Zk$ is $A$--invariant.  For the final assertion,
observe that $U(z,k)$ contains the thick set $K \subset U(z)$ posited in the definition of $\Zk$.
\qed

\bold{Thickness and moduli.}
To complete the proof Theorem \ref{thm:zk}, we just
need to relate thickness to the results of \S\ref{sec:mod}.
For the next statement, we regard $\rhat = \reals \cup \{\infty\}$ as a circle on $S^2 \isom \chat$.

\begin{prop}
\label{prop:kdelta}
Let $K \subset \Rhat$ be a Cantor set of modulus $\delta>0$ containing $\infty$.
Then $T = K \cap \reals$ is a globally $k$--thick subset of $\reals$, where $k>1$ depends only on $\delta$.
\end{prop}

\bold{Proof.}
Use Proposition \ref{prop:xratio} to relate the modulus of $K$ to the relative sizes of
gaps in $\reals-K$.
\qed

\bold{Proof of Theorem \ref{thm:zk}.}
Since $M$ is acylindrical, by Theorem \ref{thm:robust} there exists a $\delta > 0$ such that
for any $C \mem \cC^*$, there exists a Cantor set $K$ of modulus $\delta$ with
\begin{displaymath}
	K \subset C \cap \Lambda \subset S^2.
\end{displaymath}
By Proposition \ref{prop:kdelta}, there exists a $k_0$ such that 
$T \subset \reals$ is globally $k_0$--thick whenever $T \cup \infty$ is
a Cantor set of modulus $\delta$.

Let $P^*$ be a plane in $M^*$.  Choose $C \mem \cC^*$ such that
the image of $\hull(C)$ in $M^*$ contains $P^*$.
Let $K\subset C \cap \Lambda$ be the Cantor set of modulus $\delta$
provided by Theorem \ref{thm:robust}.

By a change of coordinates, we can arrange that $0,\infty \mem K \subset \Rhat$.
Let $\zt \mem \F \half^3$ be any frame tangent to $\hull(\Rhat)$ along
the geodesic $\gamma$ joining zero to infinity, and let $z$ denote its projection to $\FM$.
Then $z$ is tangent to $P^*$.
It is readily verified that there exists an isomorphism $U \isom \reals$ sending
$U(z)$ to $\reals \cap \Lambda$.   Since $0 \mem K \subset \reals \cap \Lambda$
and $K$ is globally $k_0$--thick, we have $z \mem \RF_{k_0} M$ as well.
Thus the Theorem holds for all $k \ge k_0$.
\qed

\bold{Comparison with the rigid case.}
We conclude by comparing the case of a {\em general} convex cocompact acylindrical 3-manifold $M$,
treated by Theorems \ref{thm:robust} and \ref{thm:zk},
with the {\em rigid case}, studied in \cite{McMullen:Mohammadi:Oh:sier}.

In the rigid case, every component $D_i$ of $S^2-\Lambda$ is a {\em round} disk;
hence $C \cap D_i$ is {\em connected} for all $C \mem \cC^*$, and one can show: 
\begin{quote}
	{\em $K = C \cap \Lambda$ is a compact set of definite modulus $\forall C \mem \cC^*$.}
\end{quote}
See \cite[Lemma 9.2]{McMullen:Mohammadi:Oh:sier}.
Similarly, all horocycles passing through $\RFM$ are {\em recurrent}, and
$\Zk = \RFM$ for all $k$ sufficiently large.

On the other hand, when $M$ is not rigid, there are cases where both these properties fail.
For example, suppose the bending measure of $\hull(\Lambda)$ has an atom
of mass $\theta$
along the geodesic $\gamma$ joining $p,q \mem \Lambda$.
Then we can change coordinates on $S^2 \isom \chat$ so that $p=0$, $q=\infty$,
and $\Lambda$ is contained in the wedge defined by $|\arg(z)| < \pi - \theta/2$.
Then the circle $C \mem \cC^*$ defined by $\Re(z)=1$
cannot meet the limit set in a set of positive modulus, since
$\infty$ is an isolated point of $C \cap \Lambda$.

Similarly, the horocycle in $\half^3 = \cx \times \reals_+$ defined by
$\eta(t) = (it,1)$ crosses $\gamma$ when $t=0$, and satisfies
$d(\eta(t),\hull(\Lambda)) \arrow \infty$ as $|t| \arrow \infty$.
Projecting to $M$, we obtain a {\em divergent} horocycle orbit $xU$ with $x \mem \RFM$.
In particular, $x \mem \RFM - \Zk$ for all $k$.

Nevertheless $C \cap \Lambda$ can {\em contain} a Cantor set
of positive modulus, consistent with Theorem \ref{thm:robust}.

\section{The boundary of the convex core}
\label{sec:hull}

In this short section we analyze the behavior of $C \cap \Lambda$ for circles that
meet the limit set but do not separate it.   
The result we need does not require that $M$ is acylindrical, only
that its convex core is compact.

\begin{theorem}
\label{thm:hull}
Let $M = \Gamma \bs \half^3$ be a convex cocompact 
3--manifold with limit set $\Lambda$.
Let $C$ be the boundary of a supporting hyperplane for $\hull(\Lambda)$.
Then: 
\begin{enumerate}
	\item
$\Gamma^C$ is a convex cocompact Fuchsian group; and
	\item
There is a finite set $\Lambda_0$ such that
\end{enumerate}
\begin{displaymath}
	C \cap \Lambda = \Lambda(\Gamma^C) \cup \Gamma^C \Lambda_0.
\end{displaymath}
\end{theorem}
Here $\Lambda(\Gamma^C)$ denotes the limit set of 
$\Gamma^C = \{g \mem \Gamma \st g(C) = C\}$. 

\begin{cor}
If the projection of $\hull(C)$ to $M$ gives a plane
$P$ disjoint from $M^*$ but tangent to a frame in $\Zk$,
then $\Gamma^C$ is nonelementary.
\end{cor}

\bold{Proof.}
The hypotheses guarantee that $C$ does not separate $\Lambda$,
and $C \cap \Lambda$ contains an (uncountable) Cantor set of positive modulus.
Then by the preceding result, $\Lambda(\Gamma^C)$ is uncountable,
so $\Gamma^C$ is nonelementary.
\qed

\bold{Proof of Theorem \ref{thm:hull}.}
We will use the theory of the bending lamination, 
developed in \cite{Thurston:book:GTTM},
\cite{Epstein:Marden:convex},
\cite{Kamishima:Tan} and elsewhere.

If $M^*$ is empty, then $\Lambda$ is contained 
in a circle and the result is immediate.
The desired result is also immediate if $C \cap \Lambda$ is finite,
because $\Lambda(\Gamma^C) \subset C \cap \Lambda$.

Now assume $C \cap \Lambda$ is infinite and $M^*$ is nonempty.
Then $K = \bdry \core(M)$ is a finite union of disjoint 
compact pleated surfaces with bending lamination $\beta$.
Let 
\begin{displaymath}
	f : S = \Gamma^C \bs \hull(C \cap \Lambda) \arrow M
\end{displaymath}
be the natural projection.  Since $|C \cap \Lambda| > 2$,
$S$ is a metrically complete hyperbolic surface with geodesic boundary,
with nonempty interior $S_0$.
The map $f$ sends $S_0$ isometrically to a component
of $K-\beta$; in particular, $S_0$ has finite area.    
It follows that the ends of $S_0$ consist of the regions
between finitely many pairs of geodesics which are
tangent at infinity; for an example, see Figure \ref{fig:crown}.
Consequently, we can find a finite set $\Lambda_0 \subset \Lambda$
(corresponding to the finitely many ends of $S_0$) such that
\begin{displaymath}
	C \cap \Lambda = \Lambda(\Gamma^C) \cup \Gamma^C \Lambda_0.
\end{displaymath}
The group $\Gamma^C$ is convex cocompact
because $S$ has finite area and $\Gamma$ contains no parabolic elements.
\qed

\makefig{A surface with a crown.}{fig:crown}{
	\includegraphics[height=1.2in]{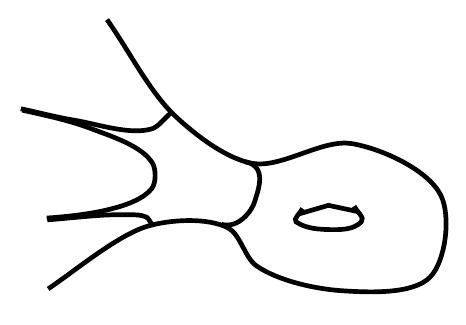}}

\section{Planes near the boundary of the convex core}
\label{sec:near}

In this section we take a step towards the proof of Theorem \ref{thm:homog}
by establishing two density results.

\begin{theorem}
\label{thm:case1}
Let $M = \Gamma \bs \half^3$ be a convex cocompact
3--manifold with incompressible boundary.
Consider a sequence of circles $C_n \arrow C$ with $C_n \mem \cC^*$ but $C \notmem \cC^*$.
Suppose that $L = \liminf (C_n \cap \Lambda)$ is uncountable.
Then ${\bigcup \Gamma C_n}$ is dense in $\cC^*$.
\end{theorem}

Under the same assumptions on $M$ we obtain:

\begin{cor}
\label{cor:case1}
Consider an $H$--invariant set $E \subset F^*$ and fix $k>1$.
If the closure of $E \cap \Zk$ 
contains a point outside $F^*$, then $E$ is dense in $F^*$.
\end{cor}

\bold{Proof.}
Consider a sequence $x_n \mem E \cap \Zk$ such that $x_n \arrow x \mem \Zk-F^*$.  We then have a corresponding
sequence of circles $C_n \mem \cC^*$ such that $C_n \arrow C \notmem \cC^*$.
(The circles are chosen so that $x_n$ is tangent to the image of $\hull(C_n)$ in $M$.)

Pass to a subsequence such that $U(x_n)$ (defined using equation (\ref{eq:Uz}))
converges, in the Hausdorff topology, to a closed set $K \subset U(x)$.
Then $C_n \cap \Lambda$ also converges, to a compact set $L \subset C$
homeomorphic to the 1-point compactification of $K$.
The fact that $x_n \mem \Zk$ implies that $K$ contains a globally $k$--thick set; hence $K$
is uncountable, so $L$ is as well.
Then by the result above, $\bigcup \Gamma C_n$ is dense in $\cC^*$, so $E$ is dense in $F^*$.
\qed

Roughly speaking, these results show that planes $P^*$ that
are nearly tangent to $\bdry M^*$ are also nearly dense in $M^*$,
subject to a condition on $\Zk$ that is automatic in the acylindrical case
by Theorem \ref{thm:zk}.

\bold{Fuchsian dynamics.}
The proof of Theorem \ref{thm:case1} exploits the dynamics of
the Fuchsian group $\Gamma^C$.
Given an open round disk $D \subset S^2$ and a closed subset
$E \subset \bdry D$, we let $\hull(E,D) \subset D$ denote the convex hull
of $E$ in the hyperbolic metric on $D$.

The principle we will use is 
\cite[Cor. 3.2]{McMullen:Mohammadi:Oh:sier},
which we restate as follows.

\begin{theorem}
\label{thm:Fuchs}
Let $M = \Gamma \bs \half^3$ be a convex cocompact hyperbolic 3--manifold.
Let $D \subset S^2$ be a round open disk that meets $\Lambda$, and let $C = \bdry D$.
Suppose $\Gamma^C$ is a nonelementary, finitely generated group, and let $C_n \arrow C$ be a sequence
of circles such that
\begin{displaymath}
	C_n \cap \hull(\Lambda(\Gamma^C),D) \neq \emptyset .
\end{displaymath}
Then the closure of $\bigcup \Gamma C_n$ in $\cC$ contains every circle that meets $\Lambda$.
\end{theorem}

\bold{Proof of Theorem \ref{thm:case1}.}
Let $D$ and $D'$ denote the two components of $S^2-C$.
Since $C \notmem \cC^*$, at least one of the components, say $D'$, is contained in $\Omega$.
Since $L \subset C \cap \Lambda$ is uncountable, $\Gamma^C$ is nonelementary
and finitely generated by Theorem \ref{thm:hull}.  Consider an ideal pentagon
\begin{equation}
\label{eq:hex}
	X = \hull(V,D) \subset \hull(\Lambda(\Gamma^C),D)
\end{equation}
whose five vertices $V$ lie in $L$.  Since $L = \liminf C_n \cap \Lambda$,
we can find `vertices'
\begin{displaymath}
	V_n \subset C_n \cap \Lambda, \;\;\;|V_n| =5,
\end{displaymath}
such that $V_n \arrow V$.
In particular, $|C_n \cap \Lambda| \ge 3$ for all $n$.

Note that $C_n$ is the unique circle passing through any three points of $V_n$.
If three of these points were to lie in $\Dbar'$, then we would have $C_n \subset \Dbar'$,
and hence $|C_n \cap \Lambda| \le 1$, since $C_n \neq C = \bdry \Dbar'$ and
$D' \subset \Omega$.  Hence $|V_n \cap D| \ge 3$.
Since $|C_n \cap C| \le 2$, at least two adjacent components of $C_n-V_n$ 
are contained in $D$.
It follows easily that $C_n$ meets $\hull(V,D)$ for all $n$ sufficiently large.
Using equation (\ref{eq:hex})
we can then apply Theorem \ref{thm:Fuchs} to conclude that 
$\bigcup \Gamma C_n$ is dense in $\cC^*$, since every $C \mem \cC^*$ meets $\Lambda$.
\qed

\section{Planes far from the boundary}
\label{sec:far}

In this section we finally prove Theorem \ref{thm:homog},
which we restate as Corollary \ref{cor:homog}.
The proof rests on the following 
more general density theorem.

\begin{theorem}  
\label{thm:case2}
Let  $M=\Gamma \backslash \bH^3$ be a convex cocompact 3-manifold with incompressible boundary.
Let $C_i \arrow C$ be a convergent sequence in $\cC^*$, with $C \notmem \bigcup \Gamma C_i$.

Suppose that $C \cap \Lambda$ contains a Cantor set of positive modulus.
Then $\bigcup \Gamma C_i$ is dense in $\cC^*$.
\end{theorem}

\begin{cor}
\label{cor:homog}
If $M  = \Gamma \bs \half^3$ is a convex cocompact acylindrical 3-manifold, then any $\Gamma$--invariant
set $E \subset \cC^*$ is either closed or dense in $\cC^*$.
\end{cor}

\bold{Proof.}
Suppose $E$ is not closed in $\cC^*$.  Then we can find a sequence $C_i \mem E$
converging to $C \mem \cC^*-E$.
Since $M$ is acylindrical, $C$ meets $\Lambda$ in a Cantor set of positive modulus, 
by Theorem \ref{thm:robust}.
Since $E$ is $\Gamma$-invariant, the preceding
result shows that $\bigcup \Gamma C_i$ is dense in $\cC^*$, so the same is true for $E$.
\qed

The proof of Theorem \ref{thm:case2}
follows the same lines as the proof of Theorem 7.3 in \cite[\S9]{McMullen:Mohammadi:Oh:sier}.
We will freely quote results from \cite{McMullen:Mohammadi:Oh:sier} in the course of the proof.
The notation from Table \ref{tab:Lie} for the subgroups $U,V,A,N$ of $G$ and other objects
will also be in play.
A generalization of Theorem \ref{thm:case2} to manifolds with compressible
boundary is stated at the end of this section.

\bold{Setup in the frame bundle.}
To prepare for the proof, we first reformulate it in terms of the frame bundle.

Let $C_i \arrow C$ as in the statement of Theorem \ref{thm:case2}.
Since $C \cap \Lambda$ contains a Cantor set of positive modulus, 
by Proposition \ref{prop:kdelta} we can choose $k>1$ and $x_\infty \mem \Zk \cap F^*$
such that $x_\infty H$ corresponds to $\Gamma C$.
Let us also choose $x_i \arrow x_\infty$ in $F^*$ such that $x_i H$ corresponds to $\Gamma C_i$.
Since $C \notmem \bigcup \Gamma C_i$, we also have
\begin{displaymath}
	x_\infty \notmem E = \bigcup x_i H.
\end{displaymath}
To prove Theorem \ref{thm:case2} we need to show:
\begin{displaymath}
	\text{\em $E$ is dense in $F^*$} .
\end{displaymath}
We may also assume that:
\begin{equation}
\label{eq:cpt}
	\text{\em The set $\Ebar \cap \Zk \cap F^*$ is compact.} 
\end{equation}
Otherwise $\Ebar \cap F^* = F^*$ by Corollary \ref{cor:case1}, and hence $E$ is dense in $F^*$.

\bold{Dynamics of semigroups.}
We say that $L \subset G$ is a {\em $1$--parameter semigroup} if there exists a nonzero $\xi \mem \Lie(G)$ such that
\begin{displaymath}
	L = \{\exp(t \xi) \st t \ge 0\} .
\end{displaymath}
To show a set is dense in $F^*$, we will use the following fact.

\begin{prop}
\label{prop:sweep}
Let $L \subset V$ be a 1--parameter semigroup.
Then $\closure{xLH}$ contains $F^*$ for all $x \mem F^*$.
\end{prop}

\bold{Proof.}
Let $C \mem \cC^*$ be a circle corresponding to $xH$.
Then $xLH$ corresponds to a family of circles $C_\alpha$
such that $\bigcup C_\alpha$ contains one of the components of $S^2-C$.
Since $C \mem \cC^*$, both components meet the limit set.
Hence $\closure{\Gamma C_\alpha} \superset \cC^*$ for some $\alpha$
by \cite[Cor. 4.2]{McMullen:Mohammadi:Oh:sier}.
\qed

\bold{The staccato horocycle flow.}
Recall that the compact set $\Zk$ is invariant under the geodesic flow $A$.
Moreover, Proposition \ref{prop:zkcompact} states that
\begin{displaymath}
        U(z,k) = \{u \mem U \st zu \mem \Zk\} 
\end{displaymath}
is a thick subset of $U$, for all $z \mem \Zk$.
In other words, $\Zk$ is also invariant under the {\em staccato horocycle flow},
which is interrupted outside of $U(z,k)$.

\bold{Recurrence.}
Next we define a compact set $W$ with
\begin{displaymath}
	x_\infty \mem W \subset \Ebar \cap F^* 
\end{displaymath}
with good recurrence properties for the horocycle flow.
Namely, we let
\begin{equation}
\label{eq:W}
	W = 
	\begin{cases}
		(\Ebar-E) \cap \Zk \cap F^* & \text{if this set is compact, and}\\
		\Ebar \cap \Zk \cap F^* & \text{otherwise} .
	\end{cases}
\end{equation}
(This definition is motivated by the proof of Lemma \ref{lem:Y}.)

In either case, $W$ is compact by assumption (\ref{eq:cpt}).
Since $\Ebar \cap F^*$ is $H$--invariant, we have
\begin{displaymath}
	WA = W \AND WU \cap \Zk \subset W .
\end{displaymath}
The second inclusion gives good recurrence; namely, we have
\begin{equation}
\label{eq:uxk}
	x U(x,k) \subset W
\end{equation}
for all $x \mem W$; and $U(x,k)$ is thick, because $W \subset \Zk$.

\bold{The horocycle flow.}
We now exploit the fact that $\Ebar$ is invariant under the horocycle flow.
The 1--parameter horocycle subgroup $U \subset H$ is distinguished by the fact that
its normalizer contains (with finite index) the 
large subgroup $AN \subset G$.
If $X$ is $U$--invariant, then so is $Xg$ for any $g \mem AN$.

\bold{Minimal sets.}
A closed set $Y$ is a {\em $U$--minimal set} for $\Ebar$ with respect to $W$ if
$Y \subset \Ebar$, $Y$ meets $W$, $YU=Y$, and
\begin{displaymath}
	\closure{yU}=Y \;\;\text{for all $y\in Y\cap W$} .
\end{displaymath}
Note that $\Ebar$ itself has all these properties except for the last.
The existence of a minimal set $Y$ follows from
the Axiom of Choice and compactness of $W$.
From now on we will assume that:
\begin{center}
	{\em $Y$ is a $U$--minimal set for $\Ebar$ with respect to $W$.}
\end{center}

To show that $\Ebar$ is large, our strategy is to 
show it contains $Yg$ for many $g \mem AN$.
To this end, we remark that for $g \mem AN$:
\begin{center}
	{\em If $(Y \cap W)g$ meets $\Ebar$, then $Yg \subset \Ebar$.}
\end{center}
Indeed, in this case by minimality we have:
\begin{equation}
\label{eq:ug}
	\Ebar \superset \closure{ygU} = \closure{yU}g = Yg,
\end{equation}
where $yg \mem (Y \cap W)g \cap \Ebar$.

\bold{Translation of $Y$ inside of $Y$.}
The fact that horocycles in $Y$ return frequently to $W$ allows one
to deduce additional invariance properties for $Y$ itself.
Note that the orbits of $AV$ are orthogonal to the orbits of $U$ in the
Riemannian metric on $\FM$.
 
\begin{lemma}
\label{lem:yw} 
There exists a $1$-parameter semigroup $L\subset AV$ such that 
\begin{displaymath}
	YL \subset Y.
\end{displaymath}
\end{lemma}

\bold{Proof.} 
In the rigid acylindrical case, this  is Theorem 9.4 in \cite{McMullen:Mohammadi:Oh:sier} for $W=\RFM$.
The only property of $\RFM$ used in the proof is 
the $k$-thickness of $\{u\in U: xu\in \RFM\}$ for any $x\in \RFM$.
Hence the proof works verbatim with $W$ replacing $\RFM$, in view of equation (\ref{eq:uxk}).  In fact $YL=Y$ since $\id \mem L$.
\qed

\bold{Translation of $Y$ inside of $\Ebar$.}
Our next goal is to find more elements $g \mem G$ 
that satisfy $Yg \subset \Ebar$.
Consider the closed set $S(Y) \subset G$ defined by
\begin{displaymath}
	S(Y) =\{g\in G: (Y\cap W) g \cap \Ebar\ne \emptyset\} .
\end{displaymath}
Since $\Ebar$ is $H$--invariant, we have $S(Y)H = S(Y)$.

\begin{lemma}
\label{lem:vn} 
If $S(Y)$ contains a sequence $g_n \to \id$ in $G-H$, then there exists $v_n\in V - \{\id\}$ tending to $\id$  such that $$Yv_n\subset \Ebar .$$
\end{lemma}

\bold{Proof.}
Let $g_n\in S(Y)$ be a sequence tending to $\id$ in $G-H$.
First suppose that there is a subsequence, which we continue to denote by $\{g_n\}$, of the form $g_n=v_nh_n\in VH$.
Since $g_n\not\in H$, we have $v_n\neq \id$ for all $n$. 
The claim then follows from the $H$-invariance of $S(Y)$ and the $U$-minimality of $Y$, see~\eqref{eq:ug}.

Therefore, assume that $g_n\notin VH$ for all large $n$. 
Since $g_n\in S(Y)$, there exist $y_n\in Y\cap W$ such that $y_n g_n \in\Ebar$. 

Since $Y$ is $U$--invariant and $WU \cap \Zk \subset \Zk$, we have
$yU(y,k) \subset Y$ for all $y \mem Y$,
and $U(y,k)$ is a $k$--thick subset of $U$.

Therefore, by \cite[Thm. 8.1]{McMullen:Mohammadi:Oh:sier}, 
for any neighborhood $G_0$ of the identity in $G$ 
we can choose $u_n \in U(y_n,k)$ and $h_n\in H$ such that 
\[
u_n^{-1} g_nh_n\to v\in V\cap G_0-\{\id\}.
\] 
After passing to a subsequence, we have $y_nu_n\to y_0\in Y\cap W$. Hence
\[
y_n g_n h_n =(y_n u_n)(u_n^{-1} g_n h_n)\in \Ebar
\] 
converges to $y_0 v\in \Ebar$.

Since $Y$ is $U$-minimal with respect to $W$ and $y_0\in Y\cap W$,
we have 
\[
\overline{y_0vU}=\overline{y_0U}v=Yv\subset \Ebar.
\]
Since $G_0$ was an arbitrary neighborhood of the identity,
the claim follows.
\qed

\bold{Choosing $Y$.}
In general there are many possibilities for the minimal set $Y$, and it may
be hard to describe a particular one, since the existence of a minimal set is
proved using the Axiom  of Choice.  The next result shows that,
nevertheless, we can choose $Y$ so it remains inside $\Ebar$ under
suitable translations transverse to $H$ but still in $AN$.

\begin{lemma}
\label{lem:Y} 
There exists a $U$--minimal set $Y$ for $\Ebar$ with respect to $W$,
and a sequence $v_n \arrow \id$ in $V-\{\id\}$, such that
\begin{displaymath}
	Y v_n \subset \Ebar
\end{displaymath}
for all $n$.
\end{lemma}

\bold{Proof.}
By Lemma \ref{lem:vn}, it suffices to show that $Y$ can be chosen so that
$S(Y)$ contains a sequence $g_n \to \id$ in $G-H$.
We break the analysis into two cases, depending on whether or not
$E$ meets the compact set $W$.

First consider the case where $E$ is disjoint from $W$.
Let $Y$ be a $U$--minimal set for $\Ebar$ with respect to $W$.
Choose $y \mem Y \cap W$.
Since $Y\subset \Ebar$, there exist $g_n \to \id$ such that $yg_n\in E$.
Then $y\notin E$, and hence $g_n \in G-H$, so we are done.

Now suppose $E$ meets $W$.  Then $W - E$ is not closed, by equation (\ref{eq:W}).
So in this case there exists a sequence
$x_n \mem W-E$ with $x_n \arrow x \mem E \cap W$.
In particular, $\overline{xH}\cap W\neq\emptyset$.
Thus there exists a $U$-minimal set $Y$ for $\overline{xH}$
with respect to $W$.

We now consider two cases.
Assume first that $Y\cap W\subset xH$.
Pick $y\in Y\cap W$; then $y=xh$ for some $h\in H$. 
Since $x_n \to x$ we have $x_n h\to y$. Now writing $yg_n=x_n h$,
we have $g_n \to \id$. As $y\in xH\subset E$ and $x_n \notin E$,
we have $g_n\in G-H$, and we are done.

Now suppose that $W\cap Y\not\subset xH$. Choose  $y \mem (W \cap Y) - xH$.
Since we have $Y\subset\overline{xH}$, 
there exist $g_n \to \id$ with $yg_n\in xH$. Moreover, $g_n \in G-H$ since $y\notin xH$,
and the proof is complete in this case as well.
\qed

\bold{Semigroups.}
We are now ready to complete the proof of Theorem \ref{thm:case2}.
We will exploit the 1-parameter semigroup $L \subset AV$ guaranteed by Lemma \ref{lem:yw}.
To discuss the possibilities for $L$, let us write the elements of $V$ and $A$ as
\begin{displaymath}
	v(s) = \Mat{1 & is \\ 0 & 1}
	\AND
	a(t) = \Mat{e^t & 0 \\ 0 & e^{-t}}.
\end{displaymath}
We then have two semigroups in $V$, defined by $V_\plusorminus = \{v(s) \st \plusorminus s \ge 0\}$,
and two similar semigroups in $A_\plusorminus$ in $A$.  It will also be useful to introduce the interval
\begin{displaymath}
	V_{[a,b]} = \{v(s) \st s \mem [a,b]\}.
\end{displaymath}
In the notation above, if $L \subset AV$ is a $1$-parameter semigroup, then either
\begin{quote}
	(i) $L = V_\plusorminus$; \\
	(ii) $L = A_\plusorminus$; or\\
	(iii) $L = v^{-1} A_\plusorminus v$, for some $v\in V$, $v \neq \id$.
\end{quote}

\bold{Proof of Theorem \ref{thm:case2}.} 
To complete the proof, it only remains to show we
have $F^* \subset \Ebar$.

Choose $Y$ and $v_n \mem V$ so that $Yv_n \subset \Ebar$
as in Lemma \ref{lem:Y}.  Write $v_n = v(s_n)$; then $s_n \arrow 0$ and $s_n \neq 0$.
Passing to a subsequence, we can assume $s_n$ has a definite sign, say $s_n>0$.

By Lemma \ref{lem:yw}, there is a $1$--parameter semigroup $L \subset AV$ such that
\begin{displaymath}
	YL\subset Y.
\end{displaymath}
The rest of the argument breaks into 3 cases, depending on whether $L$ is of
type (i), (ii) or (iii) in the list above.

\bold{(i).}
If $L=V_\plusorminus$, then we have $F^* \subset \closure{YLH} \subset \Ebar H = \Ebar$
by Proposition \ref{prop:sweep}, and we are done.

\bold{(ii).}
Now suppose $L=A_\plusorminus$.  Let 
\begin{displaymath}
	B = \{\id\} \cup \bigcup A_\plusorminus v_n A.
\end{displaymath}
Since $YL \subset Y$ and $Yv_nA \subset \Ebar A = \Ebar$ for all $n$,
we have
\begin{displaymath}
	Y B \subset \Ebar .
\end{displaymath}
Note that $a(t) v(s) a(-t) = v(e^{2t}s)$.  Consequently we have
\begin{displaymath}
	v(e^{2t}s_n) \mem B
\end{displaymath}
for all $n$ and all $t$ with $a(t) \mem L = A_\plusorminus$.

Suppose $L=A_+$.  Since $s_n \arrow 0$ and $s_n>0$,
in this case we have $V_+ \subset B$; hence $YV_+H \subset \Ebar$ and we are done as in case (i).

Now suppose $L = A_-$.  In this case at least we obtain an interval
\begin{displaymath}
	V_{[0,s_1]} \subset B.
\end{displaymath}
Choose a sequence $a_n\in A$ such that $V_+  = \bigcup a_n V_{[0,s_1]} a_n^{-1}$.
Consider $y \mem Y \cap W$.
Since $ya_n^{-1}\in W$, and $W$ is compact, after passing to a subsequence we can assume that
\begin{displaymath}
	y a_n^{-1} \arrow y_0 \mem W \subset F^* .
\end{displaymath}
We then have
\begin{displaymath}
	y_0 V_+ = \bigcup ya_n^{-1} (a_n V_{[0,s_1]} a_n^{-1}) \subset \Ebar,
\end{displaymath}
which again implies that $F^* \subset \Ebar$, by Proposition \ref{prop:sweep}.
 
\bold{(iii).}
Finally, consider the case $L = v^{-1} A_\plusorminus v$ for some $v \mem V$, $v \neq \id$. 
We then have $YB \subset \Ebar$ where
\begin{displaymath}
	B = v^{-1}A_\plusorminus v A.
\end{displaymath}
By an easy computation, $B$ contains $V_{[0,\plusorminus s]}$ for some $s>0$,
and the argument is completed as in case (ii).
\qed

\bold{The compressible case.}
In conclusion, we remark that
Theorems \ref{thm:case1} and \ref{thm:case2}
remain true without the hypothesis that
$M$ has incompressible boundary, provided we replace $\cC^*$ with
\begin{displaymath}
	\cC^\# = \{C \mem \cC^* \st \text{$C$ meets $\Lambda$}\} 
\end{displaymath}
and require that $M^*$ is nonempty.
The proofs are simple variants of those just presented.

\bibliography{math}
\bibliographystyle{math}

\bigskip
{\sc Mathematics Department, Harvard University, 1 Oxford St, Cambridge, MA 02138-2901}
\bigskip

{\sc Mathematics Department,  UC San Diego, 9500 Gilman Dr, La Jolla, CA 92093}
\bigskip

{\sc Mathematics Department, Yale University, 10 Hillhouse Avenue, New Haven, CT 06511,
	{\em and} Korea Institute for Advanced Study, Seoul,  Korea}

\end{document}